\documentclass[12pt]{amsart}
\newcommand{\R}{{\mathbf R}}
\newcommand{\e}{\epsilon}
\renewcommand{\d}{\partial}
\newcommand{\half}{\frac{1}{2}}
\newtheorem{theo}{Theorem}

\title{Gunther's proof of Nash's isometric embedding theorem}
\author{Deane Yang}
\address{Department of Mathematics\\Polytechnic University\\
Six Metrotech Center\\Brooklyn NY 11201}
\email{yang@math.poly.edu}
\begin{document}
\maketitle

\section{Preface}
Around 1987 a German mathematician named Matthias Gunther found a new
way of obtaining the existence of isometric embeddings of a Riemannian
manifold. His proof appeared in \cite{Gun89b, Gun91}.
His approach avoids the so-called Nash-Moser iteration scheme
and, therefore, the need to prove
smooth tame or Moser-type estimates for
the inverse of the linearized operator. This simplifies the proof
of Nash's isometric embedding theorem \cite{Nas56} considerably.

This is an informal expository note describing his proof. 
It was originally written, because when I first learned Gunther's proof,
it had not appeared either in preprint or published form, and I felt
that everyone should know about it. Moreover, since he
is at Leipzig, which at the time was part of East Germany, very few
mathematicians in the U.S. knew about him or his proof.

Since  many 
still seem to be unaware of Gunther's proof, even after he gave a talk
at the International Congress of Mathematicians at Kyoto in 1990
and published his proof in the proceedings \cite{Gun91}, I have
updated this note and continue to distribute it. I do, however,
encourage you to seek out Gunther's own presentations of his proof.

\section{Introduction}

Let $M$ be a smooth $n$--dimensional manifold. Given an embedding
$u: M \to \R^N$, the standard inner product on $\R^N$ induces a 
Riemannian metric on $M$. We shall denote this metric by $du\cdot du$.
In particular, given a Riemannian metric $g$ on $M$, we say that
the embedding $u$ is {\em isometric}, if
$$
du\cdot du = g
$$

Let
$N \ge \half n(n+1)$. A $C^2$ immersion $u: M \rightarrow \R^N$ is {\em
free} if for every $x \in M$, 
$$\d_iu(x), \d_i\d_ju(x),\ 1 \le i, j \le n,
$$
span a $\min(N,n + \half n(n+1))$--dimensional linear subspace of $\R^N$.

The only place where Gunther's proof differs from earlier proofs of
existence lies in showing that given a smooth, free embedding $u_0: M
\rightarrow \R^N$ and a smooth Riemannian metric $g$ sufficiently close
(in a sense to be made precise later) to $du_0\cdot du_0$, there exists a
smooth embedding $u: M \rightarrow \R^N$ close to $u_0$ such that
\begin{equation}\label{isometry}
du\cdot du = g.
\end{equation}

Although it is not necessary, we shall 
simplify the exposition by
assuming the existence of ``global'' co--ordinates on $M$. If
$M$ is compact, this is obtained by embedding $M$ smoothly into a
torus of larger dimension
 and extending smoothly the embedding $u_0$ and the metric $g$ to
the torus so that $g$ remains close to $du_0\cdot du_0$. Otherwise, if
all we are trying to prove is a local existence theorem, we can assume
that $M$ is diffeomorphic to an open set in $\R^n$. In the 
discussion below, $x^1,\ldots,x^n$ are assumed to be global
co--ordinates on $M$. (If $M$ does not have global co--ordinates,
then all the calculations below should be done using a fixed smooth
background metric $\hat{g}$, instead of the flat metric implied by
the global co--ordinates, and its Levi--Civita connection. Extra terms involving
the curvature of $\hat{g}$ and the covariant derivative of curvature
appear, but they are all of lower order and do not affect the proof at all.)

Let $v = u - u_0$ and $h = g - du_0\cdot du_0$. For convenience we
shall denote
$$
        u_i = \frac{\d u_0}{\d x^i},\ u_{ij} = \frac{\d^2u_0}{\d x^i\d
x^j}.
$$
Then (\ref{isometry})
is equivalent to
\begin{equation}\label{pert}
u_i\cdot\d_jv + u_j\cdot\d_iv + \d_iv\cdot\d_jv = h_{ij},\ 1
\le i, j \le n.
\end{equation}
Applying the standard ``integration by parts'' trick, (\ref{pert})
can be rewritten as
\begin{equation}\label{perturb}
\d_i(u_j\cdot v) + \d_j(u_i\cdot v) - 2u_{ij}\cdot v +
\d_iv\cdot\d_jv = h_{ij}.
\end{equation}
This can be written abstractly in the following form:
$$
        L_0v + Q(v,v) = h,
$$
where $L_0$ is a linear operator and $Q$ is bilinear. Nash's trick,
when $N \ge \half n(n+1) + n$,
was to observe that the linear differential operator $L_0$  could be
inverted by a zeroth order differential operator $M_0$. More recently,
M. Gromov and Bryant-Griffiths-Yang have found cases where $N < \half
n(n+1) + n$ and $L_0$ admits a right inverse $M_0$ which ``loses'' a
fixed number of derivatives. In all cases there is a loss in
regularity, so that standard implicit 
function theorems or contraction map arguments do not seem to apply.
Instead, the so--called Nash--Moser iteration scheme must be used.

Gunther's ingenious trick can be decribed as follows: He finds new
{\em nonlocal} bilinear operators $Q_1$ and $Q_2$ such that
\begin{equation}\label{divide}
        Q = L_0Q_1 + Q_2,
\end{equation}
where $Q_1$ is zeroth order and $Q_2$ is of any given negative order,
i.e. it is a bilinear smoothing operator. Actually, in the specific
situation here, the operator $Q_2$ will be identically zero.
Then the contraction mapping argument can be applied to the equation
$$
        v = M_0(h - Q_1(v,v)) - Q_2(v,v).
$$

The splitting is obtained as follows: Let
$$
        \Delta = \sum_{i=1}^n \d_i^2.
$$
Then $\Delta - 1$ is an invertible elliptic operator on $M$. Apply it
to both sides of (\ref{perturb}). Rearranging the terms and then
applying $(\Delta-1)^{-1}$ to the resulting equation, we obtain;
$$
\d_i(u_j\cdot v + Q_j(v,v)) + \d_j(u_i\cdot v + Q_i(v,v)) -
2u_{ij}\cdot v + Q_{ij}(v,v) = h_{ij},
$$
where
\begin{eqnarray*}
Q_i(v,v) &=& (\Delta-1)^{-1}(\Delta-1)v\cdot\d_iv\\
Q_{ij}(v,v) &=&  (\Delta-1)^{-1}(2\sum_{k=1}^n \d_i\d_kv\cdot\d_j\d_kv + \d_iv\cdot\d_jv -
2(\Delta-1)v\cdot\d_i\d_jv).
\end{eqnarray*}
Since $u_0$ is free, there exists a unique $\R^N$-valued bilinear
operator $Q_0$ such that $u_i\cdot Q_0 = Q_i$ and $u_{ij}\cdot Q_0 =
Q_{ij}$. The isometric embedding equation now becomes
$$
        L_0(v - Q_0(v,v)) = h,
$$
where
$$
        (L_0v)_{ij} = \d_i(u_j\cdot v) + \d_j(u_i\cdot v) -
2u_{ij}\cdot v.
$$
Given $h = h_{ij}dx^idx^j$, define $M_0h = v$, where for every $x \in
M$, $v(x)$ is the unique vector lying in the span of $u_i(x),
u_{ij}(x)$, $1 \le i, j \le n$, satisfying the following equations
\begin{eqnarray*}
u_i\cdot v &=& 0\\
-2u_{ij}\cdot v &=& h_{ij}
\end{eqnarray*}
Clearly, $M_0$ is a right inverse for $L_0$. Therefore, to solve
(\ref{perturb}), it suffices to solve the following:
$$
        v = M_0h + Q_0(v,v).
$$
Define $\Phi(v) = M_0h + Q_0(v,v)$. If $\|v\|_{2,\alpha}$, $0 < \alpha
< 1$, is sufficiently small, then $\Phi$ is a contraction mapping on a
neighborhood of $0 \in C^{2,\alpha}(M,\R^N)$.
Moreover, the linear operator $I - Q_0(v,\cdot)$ is an elliptic zeroth
order operator and therefore if $h$ is $C^{k,\alpha}$, $k \ge 2$, then
so is $v$. In particular, if $h$ is smooth, so is $v$.

We have therefore obtained the following:
\begin{theo}[Nash, Gunther \cite{Nas56,Gun89b,Gun91}]
Let $M$ be an $n$-dimensional torus and $u_0: M \rightarrow \R^N$,  $N
\ge \half n(n+1) + n$, a smooth, free immersion. Then given $0 <
\alpha < 1$, there exists $\e > 0$ (depending on $u_0$ and $\alpha$) 
 such that given any $C^{2,\alpha}$
Riemannian metric $g$, $\|g - du_0\cdot du_0\|_{2,\alpha} < \e$, there
exists a $C^{2,\alpha}$ immersion $u$ close to $u_0$ such that $du
\cdot du = g$. Moreover, if $g$ is $C^{k,\alpha}$, $2 \le k \le
\infty$, the immersion $u$
is $C^{k,\alpha}$.
\end{theo}

\bibliographystyle{amsplain}

\end{document}